\documentclass[12pt]{article}
\usepackage{amssymb}
\usepackage{epsfig}

\catcode`\@=11
\def\numberbysection{\@addtoreset{equation}{section}
        \def\theequation{\thesection.\arabic{equation}}}

\def\beq{\begin{equation}}
\def\eeq{\end{equation}}
\numberbysection
\begin{document}
\begin{titlepage}
\begin{center}
\vskip 1.in {\Large \bf On the analytic continuation of the Riemann zeta function} \vskip 0.5in P. Valtancoli
\\[.2in]
{\em Dipartimento di Fisica, Polo Scientifico Universit\'a di Firenze \\
and INFN, Sezione di Firenze (Italy)\\
Via G. Sansone 1, 50019 Sesto Fiorentino, Italy}
\end{center}
\vskip .5in
\begin{abstract}
Analyzing in detail the analytic continuation of the Riemann zeta function we
are able to generate several new identities which may be useful for applications in physics and mathematics.
\end{abstract}
\medskip
\end{titlepage}
\pagenumbering{arabic}
\section{Introduction}

The Riemann zeta function is defined as

\begin{equation}
\zeta (s) \ = \ \sum^{\infty}_{n=1} \ \frac{1}{n^s} \ \ \ \ \ \ Re \ s > 1
\label{11}
\end{equation}

The function $ \zeta (s) $ has important applications in mathematics and physics.
In mathematics, $ \zeta (s) $ is closely related to number theory,
particularly to the distribution of prime numbers \cite{1}-\cite{2}.
In physics, $ \zeta (s) $ is involved in the calculation of the partition function in Feynman path integrals.

Euler studied the function on the real numbers. Many of the values of the Riemann zeta function at even positive integers were computed by Euler. The first of them
$ \zeta (2) $ gives a solution to the Basel problem. The values at the negative
integer points are also interesting because they are rational numbers and play
an important role in the theory of modular forms.

Euler later discovered the connection between the function $ \zeta (s) $ 
and prime numbers:
\begin{equation}
\zeta (s) \ = \ \sum^{\infty}_{n=1} \ \frac{1}{n^s} \ = \prod_{p \ prime} \ \frac{1}{1-p^{-s}} \ \ \ \ \ Re \ s > 1
\label{12}
\end{equation}

Riemann extended Euler's definition to a complex variable and proved
the existence of a functional equation

\begin{equation}
\zeta (s) \ = \ \ 2^s \pi^{s-1} \ sin \left(\frac{\pi s}{2}\right) \ \Gamma(1-s) \ \zeta (1-s)
 \ \
\label{13}
\end{equation}

The zeros of the sine function imply that $ \zeta (s) $ has a simple zero
at every negative integer $ s = -2n $, known as the trivial zeros of $ \zeta (s) $.
For $ s = 1 $, the series (\ref{11}) is the diverging harmonic series. Therefore, the function
$ \zeta (s) $ is holomorphic everywhere except for a simple pole at $ s = 1 $ with
residue $ 1 $.

Then Riemann established a relationship between the nontrivial zeros and the distribution of prime numbers. The Riemann hypothesis, one of the greatest unsolved problems in mathematics, states that all nontrivial zeros lie on the
critical line ($ Re s \ = \ \frac{1}{2} $ ) \cite{3}-\cite{4}.\cite{5}.

In this article, we limit ourselves to studying in detail the analytic continuation
of the Riemann zeta function. This can be obtained through integration by
parts. By repeating the integration by parts, we can shift the validity of the
series from $ Re \ s > 1 $ to $ Re \ s > - ( p - 1 ) $. We will show that for each $ p $,
we can associate a new identity for $ \zeta (s) $.

\section{Analytical continuation of the Riemann zeta function}

The Riemann zeta function is defined by the following series

\begin{equation}
\zeta (s) \ = \ \sum^{\infty}_{n=1} \ \frac{1}{n^s} \ \ \ \ \ \ \ \ \ Re \ s > 1
\label{21}
\end{equation}

It is convergent only in the complex half-plane defined by $ Re \ s > 1 $. By analytical continuation it can be defined in the entire complex plane of the variable $ s $ using simple integrations by parts.

The first is the following

\begin{equation}
\zeta (s) \ = \ s \ \sum^{\infty}_{n=1} \ \int^{\infty}_1 \ dx \ \frac{\Theta \ ( x-n )}{x^{s+1}} \ \ \ \ \ \ \ \ \ Re \ s > 0 \label{22}
\end{equation}

This formula allows us to extend the validity of the series (\ref{21}) to the complex half-plane $ Re \ s > 0$.

If we continue with integration by parts, we arrive at the formula

\begin{equation}
\zeta (s) \ = \ s \ (s+1) \ \sum^{\infty}_{n=1} \ \int^{\infty}_1 \ dx \ \frac{
\ (x-n) \ \Theta \ ( x-n )}{x^{s+2}} \ \ \ \ \ \ \ \ \ Re \ s > -1 \label{23}
\end{equation}

In this case, the formula's validity extends to the range $ Re \ s > -1$.
If we repeat the integration by parts, we obtain the following general formula:

\begin{equation}
\zeta (s) \ = \ \frac{ s \ (s+1)... (s+p-1) }{(p-1)!} \ \sum^{\infty}_{n=1} \ \int^{\infty}_1 \ dx \ \frac{
\ (x-n)^{p-1} \ \Theta \ ( x-n )}{x^{s+p}} \ \ \ \ \ \ \ \ \ Re \ s > -(p-1) \label{24}
\end{equation}

In this way, we can extend the analytic continuation of the series (\ref{21}) to the entire complex plane of $ s $.

In this article, we will show that we can generate a new identity for $ \zeta (s) $ for each case $ p $, and we will present explicit calculations up to $ p = 12 $.

\section{Case p = 1} 

We start from the equation:

\begin{equation}
\zeta (s) \ = \ s \ \sum^{\infty}_{n=1} \ \int^{\infty}_1 \ dx \ \frac{\Theta \ ( x-n )}{x^{s+1}} \ \ \ \ \ \ \ \ \ Re \ s > 0 \label{31}
\end{equation}

To calculate the identity for $ \zeta (s) $, we subtract from the step function $ \sum^{\infty}_{n=1} \ \Theta \ ( x-n ) $ a function of x ( $ f_1 (x) $ ) in such a way as to transform the combination $ \sum^{\infty}_{n=1} \ \Theta \ ( x-n ) \ - \ f_1 (x) $ into a periodic function. In this particular case

\begin{equation}
 f_1 (x) \ = \ x \label{32}
\end{equation}

and the periodic function is equal to $ -\{x\} $, the fractional part of $ x $.

\begin{equation}
\zeta (s) \ = \ s \ \ \int^{\infty}_1 \ dx \ \frac{x}{x^{s+1}} \ - \ s \ \ \int^{\infty}_1 \ dx \ \frac{\{x\}}{x^{s+1}}\ \ \ \ \ \ \ Re \ s > 0 \label{33}
\end{equation}

The first integral can be calculated explicitly, while the second, with a
simple step, can be translated into the following form:

\begin{equation}
s \ \ \int^{\infty}_1 \ dx \ \frac{\{x\}}{x^{s+1}}\ = \ s \ \sum^{\infty}_{n=1} \ \int^1_0 \ dx \ \frac{x}{(x+n)^{s+1}}\ \ \ \ \ \ \ Re \ s > 0 \label{34}
\end{equation}

So the formula (\ref{33}) can be rewritten as

\begin{equation}
\zeta (s) \ = \  \frac{1}{s-1} \  + \ 1 \ - \ s \ \sum^{\infty}_{n=1} \ \int^1_0 \ dx \ \frac{x}{(x+n)^{s+1}}\ \ \ \ \ \ \ Re \ s > 0 \label{35}
\end{equation}

To obtain the identity for $ \zeta (s) $, we perform an integration by parts
of the second term:

\begin{eqnarray}
& \ & \sum^{\infty}_{n=1} \ \int^1_0 \ dx \ \frac{x}{(x+n)^{s+1}} \ \ = \ \frac{1}{2} \ ( \ \zeta(s+1)-1 \ ) \ + \ (s+1) \ \sum^{\infty}_{n=1} \ \int^1_0 \ dx \ \frac{x^2}{2 \ (x+n)^{s+2}} \ \nonumber \\
& \ & \ \ \label{36}
\end{eqnarray}

By repeating the integrations by parts we finally obtain the identity

\begin{equation}
\sum^{\infty}_{k=0} \ \frac{\Gamma ( s+k) }{\Gamma (k+2) \Gamma(s)} \ 
( \ \zeta(s+k)-1 \ ) \ = \ \frac{1}{s-1} \ \
 \ \ \ Re \ s > 0 \label{37}
\end{equation}

From this formula we deduce the presence of a simple pole at $ s = 1 $ for $ \zeta (s) $ with residue $1$. Special identities can be obtained,
for example for $ s = 2 $
 
\begin{eqnarray}
 ( \zeta(2) - 1 ) & + & ( \zeta(3) - 1 ) + ( \zeta(4) - 1 )
+ ( \zeta(5) - 1 ) + .... \ = \nonumber \\
& = & \sum^{\infty}_{k=0} \ \ 
( \ \zeta(k+2)-1 \ ) \ = \ 1  \label{38}
\end{eqnarray}

\section{Case p = 2}

The case $ p = 2 $ is more interesting because it is related to the general case where Bernoulli polynomials come into play

\begin{equation}
\zeta (s) \ = \ s \ (s+1) \ \sum^{\infty}_{n=1} \ \int^{\infty}_1 \ dx \ \frac{
\ (x-n) \ \Theta \ ( x-n )}{x^{s+2}} \ \ \ \ \ \ \ \ \ Re \ s > -1 \label{41}
\end{equation}

In this case, the function to be subtracted from the step function $ \sum^{\infty}_{n=1} \ (x-n) \ \Theta \ ( x-n ) $ to make it a periodic function is calculated
as follows. Take the polynomial that corresponds to the sum of the first $ n $
positive integers:

\begin{equation}
P_1 ( n ) \ = \  \sum^{n}_{i=1} \ i \ = \ \frac{n(n+1)}{2} \ \label{42}
\end{equation}

and then replace $ n \rightarrow (x-1) $

\begin{equation}
f_2(x) \ = \ \ P_1 ( n \rightarrow x-1 ) \ = \  \ \frac{x(x-1)}{2} \ \label{43}
\end{equation}

It can be shown that the function $ \sum^{\infty}_{n=1} \ (x-n) \ \Theta \ ( x-n ) \ - \ \frac{x(x-1)}{2} $ is a periodic function.

Using the same steps that solved the case $ p=1 $ we get:

\begin{equation}
\zeta (s) \ = \  \frac{1}{s-1} \  + \ \frac{1}{2} \ + \ \frac{s(s+1)}{2} \ \sum^{\infty}_{n=1} \ \int^1_0 \ dx \ \frac{x(1-x)}{(x+n)^{s+2}}\ \ \ \ \ \ \ Re \ s > -1 \label{44}
\end{equation}

With an integration by parts it can be rewritten as

\begin{eqnarray}
\zeta (s) & = &  \frac{1}{s-1} \  + \ \frac{1}{2} \ + \ \frac{s(s+1)}{12} \ ( \zeta( s+2) -1 ) \ + \ \nonumber \\
& + &  \frac{s(s+1)(s+2)}{12} \ \sum^{\infty}_{n=1} \ \int^1_0 \ dx \ \frac{3 x^2 -2 x^3}{(x+n)^{s+3}}\ \ \ \ \ \ \ Re \ s > -1 \label{45}
\end{eqnarray}

Repeating the integration by parts, we obtain the second general identity
corresponding to the case $ p = 2 $:

\begin{equation}
\zeta(s) \ = \ \frac{1}{s-1} \ + \ \frac{1}{2} \ + \frac{1}{2} \ 
\sum^{\infty}_{k=2} \ \frac{(k-1) \Gamma ( s+k) }{\Gamma (k+2) \Gamma(s)} \ 
( \ \zeta(s+k)-1 \ )  \
 \ \ \ Re \ s > -1 \label{46}
\end{equation}

Since $ Re s > -1 $ we can use this identity to compute $ \zeta(0) $

\begin{equation}
\zeta(0) \ = \ - \frac{1}{2} \ \label{47}
\end{equation}

Furthermore, one can calculate the derivative of the Riemann zeta at $ s = 0 $

\begin{equation} 
\zeta' (0) \ = \ - 1 \ + \frac{1}{2} \ \sum^{\infty}_{k=2} \ \frac{(k-1)  }{ k (k+1) } \ ( \ \zeta(k)-1 \ )\label{48}
\end{equation}

representation of $ \zeta' (0) $ in terms of $ ( \ \zeta(k)-1 \ ) $ for $ k \geq 2 $.

\section{Case p = 3 and p = 4}

Let's start with the formula

\begin{equation}
\zeta (s) \ = \ \frac{s \ (s+1) \ (s+2)}{2}  \ \sum^{\infty}_{n=1} \ \int^{\infty}_1 \ dx \ \frac{
\ (x-n)^2 \ \Theta \ ( x-n )}{x^{s+3}} \ \ \ \ \ \ \ \ \ Re \ s > -2 \label{51}
\end{equation}

The function $ f_3 (x) $ to be subtracted from the step function $ \sum^{\infty}_{n=1}
(x-n)^2 \ \Theta \ ( x-n ) $ is calculated as follows. Take the polynomial that describes the sum of the squares of the first n positive integers

\begin{equation}
P_2 ( n ) \ = \  \sum^{n}_{i=1} \ i^2 \ = \ \frac{n(n+1)(2n+1)}{6} \ \label{52}
\end{equation}

and so

\begin{equation}
f_3(x) \ = \ \ P_2 ( n \rightarrow x-1 ) \ = \  \ \frac{(x-1)x(2x-1)}{6} \ \label{53}
\end{equation}

Repeating the integrations by parts we finally obtain the following identity:

\begin{eqnarray}
\zeta (s) & = &  \frac{1}{s-1} \  + \ \frac{1}{2} \ + \ \frac{s}{12} \ - \ 
\nonumber \\
& - &  \frac{1}{12} \ \sum^{\infty}_{k=4} \ \frac{(k-2)(k-3) \Gamma( s+k)}{\Gamma(k+2) \Gamma(s)} \ ( \zeta(s+k) - 1) 
 \ \ \ \ \ \ Re \ s > -2 \label{54}
\end{eqnarray}

This identity also holds for the next series ( $ p = 4 $ ), and therefore
we can extend the validity of (\ref{54}) to $ Re \ s > -3 $.

We can check that

\begin{eqnarray}
\zeta( -2 ) & = & 0 \nonumber \\
\zeta(0) & = & - \frac{1}{2} \nonumber \\
\zeta'(0)& = & - \frac{11}{12} \ - \ \sum^{\infty}_{k=4} \ \frac{(k-2)(k-3)}{k(k+1)}
\ ( ( \zeta(k) - 1)   \ \ \ \ \ \  \label{55}
\end{eqnarray}

\section{Case p = 5 and p = 6}

Let's start with the formula

\begin{equation}
\zeta (s) \ = \ \frac{s \ (s+1) .... (s+4)}{24}  \ \sum^{\infty}_{n=1} \ \int^{\infty}_1 \ dx \ \frac{
\ (x-n)^4 \ \Theta \ ( x-n )}{x^{s+5}} \ \ \ \ \ \ \ \ \ Re \ s > -4 \label{61}
\end{equation}

The function $ f_5(x) $ to be subtracted from the step function is calculated from the following polynomial

\begin{equation}
P_4 ( n ) \ = \  \sum^{n}_{i=1} \ i^4 \ = \ \frac{1}{5} n^5 \ + \ \frac{1}{2} n^4 \ + \ \frac{1}{3} n^3 - \frac{1}{30} n \label{62}
\end{equation}

Therefore

\begin{equation}
f_5 ( x ) \ = \ \ P_4 ( n \rightarrow x-1 ) \ = \  \ 
\frac{1}{30} ( 6 x^5 -15 x^4 +10 x^3 -x ) \ \label{63}
\end{equation}

From integration by parts we obtain the following identity

\begin{eqnarray}
\zeta (s) & = &  \frac{1}{s-1} \  + \ \frac{1}{2} \ + \ \frac{29}{360} s \ - \ 
\frac{s^2}{240} \ - \ \frac{s^3}{720} \ + \
\nonumber \\
& + &  \frac{1}{720} \ \sum^{\infty}_{k=6} \ \frac{(k-2)(k-4)(k-5)(k+9) \Gamma( s+k)}{\Gamma(k+2) \Gamma(s)} \ ( \zeta(s+k) - 1) 
 \ \ \ \ \ \ Re \ s > -4 \nonumber\\
 & \ &  \label{64}
\end{eqnarray}

In particular

\begin{equation}
\zeta(2) \ = \ \frac{49}{30} \ + \ \frac{1}{720} \ \sum^{\infty}_{k=6} (k-2)(k-4)(k-5)(k+9) \ ( \zeta(k+2) - 1)
\ \label{65}
\end{equation}

The case $ p = 6 $ leads to the same identity (\ref{64}), which therefore also holds
for $ Re \ s > - 5 $. It can be checked that this formula contains the trivial zeros for $ s = -2,-4 $.

\section{Case p = 7 and p = 8}

Let's start with the formula

\begin{equation}
\zeta (s) \ = \ \frac{s \ (s+1) .... (s+6)}{6!}  \ \sum^{\infty}_{n=1} \ \int^{\infty}_1 \ dx \ \frac{
\ (x-n)^6 \ \Theta \ ( x-n )}{x^{s+7}} \ \ \ \ \ \ \ \ \ Re \ s > -6 \label{71}
\end{equation}
The corresponding polynomial is the following

\begin{equation}
P_6 ( n ) \ = \  \sum^{n}_{i=1} \ i^6 \ = \ \frac{1}{7} n^7 \ + \ \frac{1}{2} n^6 \ + \ \frac{1}{2} n^5 - \frac{1}{6} n^3 \ + \ \frac{1}{42} n \label{72}
\end{equation}

Therefore

\begin{equation}
f_7 ( x ) \ = \ \ P_6 ( n \rightarrow x-1 ) \ = \  \ 
\frac{1}{7} x^7 \ - \ \frac{1}{2} x^6 \ + \ \frac{1}{2} x^5 \ - \ \frac{1}{6} x^3
\ + \ \frac{x}{42} \ \label{73}
\end{equation}

The solution with $ p = 7 $ is

\begin{eqnarray}
& \ & \zeta (s) \ = \  \frac{1}{s-1} \  + \ \frac{15120 + 2460 s -76 s^2 -7 s^3 + 10 s^4 + s^5 }{30240} \ - \ \nonumber \\
& - &  \frac{1}{6! 42} \ \sum^{\infty}_{k=8} \ \frac{(k-2)(k-4)(k-6)(k-7)(k^2+10k+45) \Gamma( s+k)}{\Gamma(k+2) \Gamma(s)} \ ( \zeta(s+k) - 1) 
 \ \ \ \ \ \  \nonumber\\
 & \ & Re \ s > -6  \label{74}
\end{eqnarray}

The case $ p = 8 $ leads to the same identity (\ref{74}), which therefore also holds for $ Re \ s > -7 $. It can be checked that this formula contains the trivial zeros for $ s = -2, -4, -6 $.

\section{Case p = 9 and p = 10}

Let's start with the formula

\begin{equation}
\zeta (s) \ = \ \frac{s \ (s+1) .... (s+8)}{8!}  \ \sum^{\infty}_{n=1} \ \int^{\infty}_1 \ dx \ \frac{
\ (x-n)^8 \ \Theta \ ( x-n )}{x^{s+9}} \ \ \ \ \ \ \ \ \ Re \ s > -8 \label{81}
\end{equation}

The corresponding polynomial is the following

\begin{equation}
P_8 ( n ) \ = \  \sum^{n}_{i=1} \ i^8 \ = \ \frac{1}{9} n^9 \ + \ \frac{1}{2} n^8 \ + \ \frac{2}{3} n^7 - \frac{7}{15} n^5 \ + \ \frac{2}{9} n^3 \ - \ \frac{1}{30} n \label{82}
\end{equation}

Therefore

\begin{equation}
f_9 ( x ) \ = \ \ P_8 ( n \rightarrow x-1 ) \ = \  \ 
\frac{1}{90} \ x \ ( 10 x^8 -45 x^7 + 60 x^6 -42 x^4 + 20 x^2 - 3 \ ) \ \label{83}
\end{equation}

This case is complex, and we used the MATHEMATICA program to solve it. Computer calculations show that the solution with $ p = 9$ is an eighth-degree polynomial, which fortunately can be decomposed as follows:

\begin{eqnarray}
\zeta (s) & = &  \frac{1}{s-1} \  + \ \frac{604800 + 97680 s - 4804 s^2 - 1904 s^3 - 335 s^4 - 135 s^5 - 21 s^6 - s^7 }{1209600} \ + \ \nonumber \\
& + &  \frac{1}{8! 30} \ \sum^{\infty}_{k=10} \ (k-2)(k-4)(k-6)(k-8)(k-9)
\ \cdot \nonumber \\
& \cdot & \frac{(k+5)(k^2+4k+35) \Gamma( s+k)}{\Gamma(k+2) \Gamma(s)} \ ( \zeta(s+k) - 1) 
 \ \nonumber \\
 & \ &  Re \ s > -8 \label{84}
\end{eqnarray}

The case $ p = 10 $ leads to the same identity (\ref{84}), which therefore also holds for $ Re \ s > -9 $. It can be checked that this formula contains the trivial zeros for $ s = -2, -4, -6, -8 $.

\section{Case p = 11 and p = 12}

Let's start with the formula

\begin{equation}
\zeta (s) \ = \ \frac{s \ (s+1) .... (s+10)}{10!}  \ \sum^{\infty}_{n=1} \ \int^{\infty}_1 \ dx \ \frac{
\ (x-n)^{10} \ \Theta \ ( x-n )}{x^{s+11}} \ \ \ \ \ \ \ \ \ Re \ s > -10 \label{91}
\end{equation}

The corresponding polynomial is the following

\begin{equation}
P_{10} ( n ) \ = \  \sum^{n}_{i=1} \ i^{10} \ = \ \frac{1}{11} n^{11} \ + \ \frac{1}{2} n^{10} \ + \ \frac{5}{6} n^9 - n^7 \ + \  n^5 \ - \ \frac{1}{2} n^3 \ + \ \frac{5}{66} n \label{92}
\end{equation}

Therefore

\begin{equation}
f_{11} ( x ) \ = \ \ P_{10} ( n \rightarrow x-1 ) \ = \  \ 
\frac{1}{66} \ x \ ( 6 x^{10} -33 x^9 +  55 x^8 -66 x^6 + 66 x^4 - 33 x^2 + 5 \ ) \ \label{93}
\end{equation}

This case is particularly complex, and we used the computer to solve it. In this case, it turns out that the solution with $ p = 11$ is a tenth-degree polynomial, which fortunately can be decomposed into smaller polynomials:

\begin{eqnarray}
\zeta (s) & = & \frac{1}{s-1} \  + \ \frac{119750400 + 19542240 s - 403272 s^2 + 213628 s^3 + 270090 s^4}{239500800} \ + \ \nonumber \\
& + & \frac{ 85515 s^5 + 18522 s^6 + 2532 s^7 + 180 s^8 + 5 s^9
}{239500800} \ - \ \nonumber \\
& - &  \frac{1}{10! 66} \ \sum^{\infty}_{k=12} \ (k-2)(k-4)(k-6)(k-8)(k-10)(k-11)
\ \cdot \nonumber \\
& \cdot & \frac{(5k^4 +30k^3+232k^2+1122k+2835) \Gamma( s+k)}{\Gamma(k+2) \Gamma(s)} \ ( \zeta(s+k) - 1) 
 \ \nonumber \\
 & \ &  Re \ s > -10 \label{94}
\end{eqnarray}

The case $ p = 12 $ leads to the same identity (\ref{94}), which therefore also holds
for $ Re \ s > -11 $. We can verify that this formula contains the trivial zeros for
 $ s = -2,...,-10 $. We also checked the values for $ s $ negative integer, within
 the validity range, and the formula works.

\section{Conclusions}

In this article, we studied in detail the analytic continuation of the Riemann zeta function. Analytic continuation can be performed by applying integration by parts to the definition formula. There are several ways to represent it, labeled by a number $ p $, and each way gives rise to a different identity for $ \zeta (s) $. All these ways include step functions. To calculate the corresponding identity, a Bernoulli polynomial is subtracted from the step function to convert it to a periodic
function. The periodic function then gives rise, by integration by parts to a series
that include the combinations $ ( \ \zeta(k) -1 \ )$ for $ k \geq 2 $. These values
are small and tend to zero in th limit $ k \rightarrow \infty $. The simplest identity

\begin{equation} \sum^\infty_{k=2} \ ( \ \zeta(k) -1 \ ) \ = \ 1 \label{101} 
\end{equation}

can be easily verified on a computer. As $ p $ increases, the identity for
$ \zeta (s) $ becomes increasingly complex, but in any case it is always possible
to obtain a general formula. It follows that we can obtain a representation
of the derivative of $ \zeta (s) $ as a summation over the combinations $ ( \ \zeta(k) -1 \ )$ for $ k \geq 2 $. We hope these identities can shed light on new properties of $ \zeta (s) $.

\end{document}